\date{February 18, 2010}

\documentclass[12pt]{article}
\usepackage{amsmath,amsthm,amsfonts,amssymb}   
\def\A{{\mathcal A}}
\def\B{{\mathcal B}}
\def\C{{\mathcal C}}
\def\M{{\mathcal M}}
\def\e{{\varepsilon}}
\def\Z{{\mathbb Z}}

\newtheorem{theorem}{Theorem}[section]

\newtheorem{corollary}[theorem]{Corollary}

\newtheorem{remark}[theorem]{Remark}

\title{On the Jones Index Values \\
for Conformal Subnets}

\author{
{\sc Sebastiano Carpi}\footnote{Supported in part by PRIN-MIUR, GNAMPA-INDAM and EU network ``Noncommutative Geometry" 
MRTN-CT-2006-0031962}\\  
Dipartimento di Scienze,
Universit\`a di Chieti-Pescara ``G. d'Annunzio''\\
Viale Pindaro, 42, I-65127 Pescara, Italy\\
E-mail: {\tt carpi@sci.unich.it}\\
\phantom{X}\\{\sc Yasuyuki Kawahigashi}\footnote{Supported in part by the
Grants-in-Aid for Scientific Research, JSPS.}\\
Department of Mathematical Sciences\\
University of Tokyo, Komaba, Tokyo, 153-8914, Japan\\
E-mail: {\tt yasuyuki@ms.u-tokyo.ac.jp}\\
\phantom{X}\\
{\sc Roberto Longo}$^*$
\\
Dipartimento di Matematica,
Universit\`a di Roma ``Tor Vergata'',\\
Via della Ricerca Scientifica, 1, I-00133 Roma, Italy\\
E-mail: {\tt longo@mat.uniroma2.it}}

\begin{document}

\maketitle

\begin{abstract} We consider the smallest values taken by the Jones index for an
inclusion of local conformal nets of von Neumann algebras on $S^1$ and show that
these values are quite more restricted than for an arbitrary inclusion of factors.
Below $4$, the only non-integer admissible value is $4\cos^2 \pi/10$, which is known
to be attained by a certain coset model. Then no index value is possible in the
interval between $4$ and $3 +\sqrt{3}$. The proof of this result based on
$\alpha$-induction arguments. In the case of values below 4 we also give a second
proof of the result. In the course of the latter proof we classify all possible
unitary braiding symmetries on the $A\,D\,E$ tensor categories, namely the ones
associated with the even vertices of the $A_n$, $D_{2n}$, $E_6$, $E_8$ Dynkin
diagrams.

\end{abstract} 

\thanks{\footnotesize{Supported in part by the ERC Advanced Grant 227458  
OACFT ``Operator Algebras and Conformal Field Theory"}}

\section{Introduction}
\label{Intro}

A celebrated theorem by V.F.R. Jones \cite{J} states that for an inclusion of factors
$N\subset M$ the possible values for the index are $[M:N] = 4\cos^2 \pi/n$, $n=3,
4,\dots$ or $[M:N]\geq 4$. All the values $4\cos^2 \pi/n$ in the discrete series are
realized and the in this case the inclusion is automatically irreducible, i.e.
$N'\cap M=\mathbb C$. It is also rather easy to construct a reducible inclusion of
factors with any index value larger or equal to $4$.  All the values above $4$ are
realized for non-hyperfinite irreducible subfactors as in \cite{P0}. Whether every
value larger than $4$ can be realized as the index of an irreducible inclusion of
hyperfinite factors is still an open problem, we mention that deep progresses have
been performed in the classification of finite depth subfactors, in particular by
Popa, Haagerup and Asaeda, see \cite{P,H,AH,AY}.

In this paper we begin to study the possible values of the index for an inclusion
$\A\subset \B$ of local conformal subnet on $S^1$. Indeed only M\"obius covariance
will enter in our discussion. Note that here, if the index of $\A\subset \B$ is
finite, then the inclusion $\A\subset \B$ is automatically irreducible, i.e.
$\A(I)'\cap\B(I)=\mathbb C$ for any interval $I$ of $S^1$, see \cite[Corollary 
2.7]{DLR}.

We shall see that the possible index values for $\A\subset \B$ below $3 +\sqrt{3}$
are $1,2, 3$, $4\cos^2 \pi/10$ and $4$. So only one non-integer value in the Jones
discrete series appears here: the value $4\cos^2 \pi/10$. An inclusion of conformal
net with index $4\cos^2 \pi/10$ has been considered by Rehren \cite{R}. This is
pointed out in \cite{D} where one can also find further partial results related to
this paper.

The integer index values are easily realized as orbifold, indeed if $\A$ is the
fixed-point net of $\B$ under a finite gauge symmetry group $G$ then the index of
$\A\subset \B$ is the order of $G$, so it suffices to consider for example the cyclic
orbifold of $\C\otimes\cdots\otimes\C$ for any conformal net $\C$.

One might compare our result with the restriction on the index values for a sector of
a local conformal net \cite{L0,R1}. Also in this case there are definitely more index
value gaps than in the subfactor case (in particular only four index values between
$4$ and $6$ are possible), yet there is no obvious direct link between the two cases.

We shall give two proofs of our results. The first one is based on $\alpha$-induction
arguments, has the advantage of being shorter and works for all the index values
below $3 +\sqrt{3}$.  The second proof is, in a sense, more conceptual but it works
only for the index values below $4$. It can be viewed as a result for braided tensor
categories. The basic idea is the following. Since local extensions of a conformal
net correspond to local $Q$-systems in the tensor category of DHR sectors of the net
\cite{L,LR}, our problem is of purely tensor categorical nature and it is sufficient
to consider any local conformal net whose representation tensor category contains our
given braided tensor category. We basically shall use $SU(2)_n$-models to this end.

Indeed we shall need to classify all possible unitary braiding symmetries for the
tensor categories associated with the even vertices of a Dynkin diagram of type
$A_n$, $D_{2n}$, $E_6$ and $E_8$ (the Dynkin diagrams associated with a subfactor
with index less than $4$). This result, which has its own interest, might be known to
some extent; it can be derived by means of Izumi quantum double analysis \cite{I},
but is apparently not available in the literature, cf. however \cite{KW}.

\section{$\alpha$-induction and restriction on the index values}
\label{alpha}

We consider possible index values of inclusions of
local conformal nets.
Let $\A(I)\subset \B(I)$ be an inclusion
of local conformal nets with finite index.
(See \cite{EK2} for basic notions and results
of the Jones index theory \cite{J}.)
Here the M\"obius covariance 
for $\A, \B$ is enough and we do not need
${\mathrm{Diff}}(S^1)$-covariance here.
We use the machinery of $\alpha$-induction, which has been
introduced in \cite{LR} and studied in \cite{X1,BE}.
(Strong additivity of a net is sometimes assumed, but we
do not need this assumption.)

Consider the principal and dual principal graphs of the
subfactor $\A(I)\subset \B(I)$ for a fixed interval
$I\subset S^1$.
We label $\A$-$\A$, $\B$-$\A$, and $\B$-$\B$ sectors by
symbols $\lambda$, $b$, $\tau$, respectively.  Note that
the $\A$-$\A$ sectors arising from this subfactor
are DHR sectors of the local net $\A$.  We
denote the inclusion map of this subfactor by $\iota$
and regard it as a $\B$-$\A$ sector.  The dual canonical
endomorphism $\theta$ is given by $\bar \iota \iota$
and it is regarded as an $\A$-$\A$ sector.

Fix an irreducible $\A$-$\A$ sector $\lambda$ (until the
end of this paragraph).
Let $\iota\lambda=\sum_b n_b b$ be the irreducible decomposition
of the $\B$-$\A$ sector $\iota\lambda$, where $n_b$ is a nonnegative
integer representing the multiplicity.
We also let $\alpha_\lambda=\sum_\tau m_\tau \tau$
and $b\bar\iota=\sum_\tau k_{b,\tau}\tau$ be the irreducible
decompositions, where $b$ is an arbitrary $\B$-$\A$ sector.
Note that we have
$$n_b=\langle \iota\lambda, b \rangle=
\langle \alpha_\lambda \iota, b \rangle=
\langle \alpha_\lambda, b \bar\iota\rangle=
\sum_\tau k_{b,\tau}m_\tau.$$
By locality, the $\alpha$-induction has the following
property.
$$ \langle \iota \lambda, \iota \lambda\rangle
=\langle \alpha_\lambda, \alpha_\lambda\rangle,$$
as in \cite[Proof of Theorem 3.3]{X1}. 
(Also see \cite[Part I, Theorem 3.9]{BE}.)
This implies the following.
$$\sum_b(\sum_\tau k_{b,\tau}m_\tau)^2=
\sum_\tau m_\tau^2.$$
By expanding the left hand side, we have
$$\sum_{b,\tau,\tau'}
k_{b,\tau} k_{b,\tau'} m_\tau m_{\tau'}=
\sum_\tau m_\tau^2,$$
which implies
\begin{equation}
\sum_\tau(\sum_b k_{b,\tau}^2)m_\tau^2\le
\sum_\tau m_\tau^2.
\label{xxx}
\end{equation}
For any $\tau$ with $m_\tau > 0$, there exists $b$
with $b\prec \tau \iota$, and for this $b$, we have
$k_{b,\tau}=\langle \tau, b\bar\iota\rangle=
\langle \tau \iota, b \rangle > 0$.
Then the above inequality (\ref{xxx}) implies
that there exists only one $b$ with $k_{b,\tau} > 0$, and for
this $b$, we have $k_{b,\tau}=1$.  This shows that if
we have $m_\tau > 0$, then there exists a unique irreducible $b$
with $\tau\iota=b$.

Now choose an arbitrary irreducible $\B$-$\A$ sector $b$
appearing in the dual principal graph.
Then there exists $\lambda$ with
$b\prec \iota\lambda=\alpha_\lambda\iota$.
Since $$1\le \langle \iota\lambda, b\rangle=
\langle \alpha_\lambda\iota, b\rangle
=\langle \alpha_\lambda, b\bar\iota\rangle,$$
there exists an irreducible $\B$-$\B$ sector $\tau$ with
$\tau\prec\alpha_\lambda$ and $\tau\prec b\bar\iota$.
By the arguments in the previous paragraph, we have
$\tau\iota=b$.  This has proved the following.

\begin{theorem}\label{rest}
In the above setting,
for any odd vertex $b$ of the dual principal graph
of the subfactor $\A(I)\subset \B(I)$, we have an even
vertex $\tau$ of the dual principal graph that is connected
only to $b$.
\end{theorem}

This gives the following corollary immediately, since if 
we take the
odd vertex having distance 3 from the initial vertex $*$
as $b$ in the above theorem, then the property of the
theorem is violated.

\begin{corollary}
\label{graph}
Consider a bipartite graph $G$ with the initial vertex $*$.
Suppose $G$ has a vertex with valency larger than 2 and
let $d$ be the distance from $*$ to the nearest
vertex with valency larger than 2.  If $d > 3$, then
the graph $G$ is not the dual principal graph of a
subfactor arising from an irreducible local extension
of a local conformal net.
\end{corollary}

If the Jones index is less than 4, the principal graph must
be one of the $A$-$D_{2n}$-$E_{6,8}$ diagrams.  Also in this case,
the principal graph and the dual principal graph are the same.
Only $A_2$, $A_3$, $A_5$, $D_4$ and $D_6$ satisfy the condition
in Theorem \ref{rest}.  Note that for any finite groups $H\subset G$,
we can realize $G$ as a subgroup of the automorphism group of a
local conformal net $\M$ through permutations of the tensor components.
(That is, we embed $G$ into some symmetric group $S_n$ of order $n$ and
realized $\M$ as the $n$th tensor power of some local conformal net.
Then $G$ acts on $\M$ through permutations.)  Then
we have an inclusion of local conformal net
$\A=\M^G\subset \M^H=\B$.  The principal graphs
$A_2$, $A_3$, $A_5$, $D_4$ can be realized with the choices
of $(G,H)=(\{e\}, \{e\}), (S_2, \{e\}), (S_3, S_2),
({\mathbb Z}/3{\mathbb Z}, \{e\})$, where $e$ denotes the identity
element in a group.

We know that the coset construction for the diagonal embedding
$SU(2)_3\subset SU(2)_2\otimes SU(2)_1$ gives the Virasoro net
${\mathrm{Vir}}_{7/10}$.  (See \cite[Section 3]{KL} based
on \cite{X2}.)  Then the embedding 
$SU(2)_3\otimes {\mathrm{Vir}}_{7/10}\subset SU(2)_2\otimes SU(2)_1$
has principal graph $D_6$ as noted in
\cite[(5.6)]{R}.  Thus the principal
graph $D_6$, with the index value $4\cos^2 \pi/10$, is realized for
an inclusion of local conformal nets.

Next we consider the case beyond index 4.
Haagerup \cite{H} has considered possible principal graphs for the
index range $(4,3+\sqrt3)$ and shown that such graphs are
severely restricted.  Today it has been determined that three
of his pairs of finite graphs are realized and the others are not.
(See \cite{AH,AY,BMPS} and references there.)
It is easy to see that none of Haagerup's
graphs survive the restrictions given in Corollary \ref{graph}.
(Note that Corollary \ref{graph} applies also to an infinite
graph and excludes $A_\infty$, which is in the list of
Haagerup \cite{H}.)

At the index value $3+\sqrt3$, we have an inclusion of local
conformal nets arising from conformal embedding
$SU(2)_{10}\subset SO(5)_1$, which was studied in \cite{X1}.
(Note that this subfactor is isomorphic to the well-known
GHJ-subfactor \cite[Section 4.5]{GHJ} arising from $E_6$, up
to tensoring the common injective type III$_1$ factor, as
shown in \cite[Proposition A.3]{BEK2}.)

By these considerations, we obtain the following theorem.

\begin{theorem}
\label{rest-ind}
The smallest five values of the Jones indices of inclusions
of local conformal nets are $1$, $2$, $3$,
$4\cos^2 \pi/10$, $4$ and $3+\sqrt3$.
\end{theorem}

\section{Classification of braidings}

We now present a different method to determine which index values
below 4 are possible for inclusions of local conformal nets.
Note that if an inclusion $\A(I)\subset \B(I)$ of local conformal
nets has index below 4, the system of $\A(I)$-$\A(I)$ bimodules
given by the even vertices of the principal graph has to have
a braiding, since this gives a full subcategory of the DHR-category
of the representations of the local conformal net $\A(I)$.
The principal graph in this case must be one of the 
$A_n$, $D_{2n}$, $E_6$, $E_8$ diagrams, so we will classify
all possible braidings on the even vertices of these graphs
in this section.

Ocneanu realized that one can classify braidings by studying
the ``quantum double'' system \cite{O}, and we follow the
formulation of Izumi in \cite{I}, based on the Longo-Rehren
subfactor \cite{LR}.
(We actually follow a formulation in \cite{BEK3}, where
a definition dual to that in \cite{I} is used.)
Consider a system $\Delta$ of irreducible
endomorphisms with finite indices for a type III factor $M$.
If we have a braiding $\e(\lambda,\mu)$ on the system $\Delta$,
it also gives a half-braiding naturally, so we have
a system $\{(\lambda,\e(\lambda,\cdot))\mid\lambda\in\Delta\}$
of irreducible objects,
which gives a subcategory of the quantum double category,
and this system has the same fusion rules as $\Delta$
by \cite{I}  Conversely, if we have a subcategory of the quantum
double category generated by a system
$\{(\lambda,{\mathcal E}_\lambda(\cdot))\mid\lambda\in\Delta\}$
of irreducible objects,
where ${\mathcal E}_\lambda(\cdot)$ is a half-braiding, and
this system has the same fusion rules as $\Delta$, then
it gives a braiding on $\Delta$ by
\cite[Theorem 4.6 (ii), (iii)]{I}.  Hence, for computing
the number of braidings, up to equivalence,
it is enough to count the number
of embeddings of $\Delta$ to 
$\{(\lambda,{\mathcal E}_\lambda(\cdot))\mid\lambda\in\Delta\}$,
with the same fusion rules in the quantum double category.
The quantum double category for the $A_n$, $D_{2n}$, $E_6$, $E_8$
subfactors have been computed in \cite{O}, \cite{EK1}, \cite{I},
\cite{BEK3}, so we simply count the number of such embeddings below.

We deal with five cases separately, and the first three cases
are further divided into subcases.

{\bf Case I}: $A_{2n}$.

{\bf Subcase Ia}: $A_2$.

In this case, the subfactor has index $1$, and 
we have only one even vertex.  So the number of braidings is
trivially $1$.

{\bf Subcase Ib}: $A_{2n}$ with $n>1$.
It is known that the braiding arising from $SU(2)_{2n-1}$ restricted
on the even vertices is nondegenerate.  (See \cite{O},
\cite{EK1}, \cite{I}.)  So the quantum double system is simply
given by doubling as in \cite[Corollary 7.2]{I}, and the
irreducible objects are labeled with pairs $(j,k)$ with
$j,k\in\{0,2,4,\dots,2n-2\}$ where the irreducible objects of
the $SU(2)_{2n-1}$ category are labeled with $0,1,2,\dots,2n-1$
as usual.  Then the odd vertices of the dual principal graph
of the Longo-Rehren subfactor $M\otimes M^{\mathrm{opp}}
\subset R$ are labeled with $l\in\{0,2,4,\dots,2n-2\}$, and
the number of the edges connecting $(j,k)$ and $l$ is given
by the structure constant $N_{jk}^l$ by \cite[Section 4]{I}.
Then it is easy to see that only possible embeddings are
given by $j\mapsto (j,0)$ and $j\mapsto (0,j)$.  So the
number of the braidings is 2.

{\bf Case II}: $A_{2n+1}$.  We again label the irreducible
objects of the $SU(2)_{2n}$ category with $0,1,2,\dots,2n$.
The braiding arising from $SU(2)_{2n}$ restricted
on the even vertices is degenerate and the quantum double
system is explicitly described, as in 
\cite{EK1}, \cite[Section 7]{I}.

{\bf Subcase IIa}: $A_3$.
The even vertices of $A_3$ are given by the group $\Z/2\Z$.
We draw the induction-restriction graph for the quantum double
system for the Longo-Rehren subfactor $M\otimes M^{\mathrm{opp}}
\subset R$.  Then it is easy to see that we have 2 embeddings.

{\bf Subcase IIb}: $A_5$.  It was found by Ocneanu that the
number of braidings is 3.  We present arguments here for the
sake of completeness as follows.

The dual principal graph of the Longo-Rehren
subfactor $M\otimes M^{\mathrm{opp}} \subset R$ for this case
is given in \cite[Figure 21]{EK1}.  It is easy to see that the vertex
labeled with $2$ can be mapped to $3$ irreducible sectors 
labeled with $02$, $20$, $22_-$ and
each indeed gives an embedding in the above sense.
Thus the number of the
braidings is $3$.  The restriction of the original braiding
on $SU(2)_4$ gives 2 of them, and the other comes from 
realization of the even vertices of $A_5$ from the dual
of $S_3$, which gives a completely degenerate braiding,
where all monodromy operators are trivial.

{\bf Subcase IIc}: $A_{2n+1}$ with $n>2$.
We label the irreducible objects of the quantum double system
with $(j,k)$ with $j,k\in \{0,1,2,\dots,2n\}$ satisfying
$j+k\in 2\Z$ and $(j,k)\neq(n,n)$ with identification
$(j,k)=(2n-j,2n-k)$, and $(n,n)_+$, $(n,n)_-$
as in \cite[Section 4]{EK1}, \cite[Section 7]{I}.
The statistical dimension of $(j,k)$ is
$$\sin((j+1)\pi/(2n+2))\sin((k+1)\pi/(2n+2))/\sin^2(\pi/(2n+2))$$
and that of $(n,n)_\pm$ is $1/2\sin^2(\pi/(2n+2))$, and now $n>2$, so
only $(0,2)$, $(0,2n-2)$, $(2,0)$, $(2n-2, 0)$ have the
same statistical dimensions as the irreducible object $2$
of $SU(2)_{2n}$.  From this, we see that
the only possible embeddings are
given by $j\mapsto (j,0)$ and $j\mapsto (0,j)$.  So the
number of braidings is 2.

{\bf Case III}: $D_{2n}$.  Note that the even vertices of
the $D_{2n}$ diagram are realized as the DHR sectors of
the local extensions of the $SU(2)_{4n-4}$-nets with index 2
as in \cite[Part II]{BE}.  This gives at least
2 non-degenerate embeddings.

{\bf Subcase IIIa}: $D_4$. 
The even vertices of $D_4$ are given by the group $\Z/3\Z$.
By similar arguments to the case $A_3$, we see we have 3
embeddings.  Two of them come from the local extension of
the $SU(2)_4$-net, and the other is a degenerate one coming
from $\Z/3\Z\cong\widehat{\Z/3\Z}$.

{\bf Subcase IIIb}: $D_6$.
The even vertices of $D_6$ are realized as the
irreducible DHR sectors of a local extension of $SU(2)_8$
of index 2 with the dual canonical endomorphism
$\theta=0\oplus 8$.  The $\alpha$-inductions of 0 and 2
are irreducible, and that of 4 has an irreducible 
decomposition into two pieces, so we label them with
$0,2,4_+,4_-$, respectively.  Then the fusion rules are 
commutative and non-trivial ones are as follows.
\begin{eqnarray*}
2\cdot 2&=&0+2+4_++4_-,\\
2\cdot 4_+&=& 2+4_+,\\
2\cdot 4_-&=& 2+4_-,\\
4_\pm\cdot 4_\pm&=& 0+4_\pm,\\
4_\pm\cdot 4_\mp&=& 2.
\end{eqnarray*}

Then the even vertices of the dual principal graph of the
Longo-Rehren subfactor $M\otimes M^{\mathrm{opp}}\subset R$
are labeled with pairs $(j,k)$ with $j,k\in\{0,2,4_+, 4_-\}$.
The odd vertices of the dual principal graph are labeled
with $l\in\{0,2,4_+, 4_-\}$ and the number of edges connecting
$(j,k)$ and $l$ is given by the structure constant $N_{jk}^l$.
We now count the number of
embeddings of the system $\{0,2,4_+, 4_-\}$
into the quantum double system.
From the above fusion rules, it is easy to see that 
$2$ has to be mapped to one of
$(0,2), (2,0), (4_+, 4_-), (4_-, 4_+)$.  Then it is also
easy to see that all of these give embeddings, so the
number of the braidings is 4.

These four braidings are interpreted as follows.  Since
the braiding of the even vertices of $A_4$ is nondegenerate
as see above, its quantum double system is simply a self-doubling.
That is, if we label the even vertices of $A_4$ with $1$, $\sigma$,
where $1$ is the identity sector, then the even vertices of
$D_6$ are labeled as $1\otimes 1$, $1\otimes\sigma$, $\sigma\otimes1$,
and $\sigma\otimes\sigma$.  Since the system $\{1,\sigma\}$ has
braidings $\e^\pm$ arising from $SU(2)_3$, the system
$\{1\otimes 1, 1\otimes\sigma, \sigma\otimes1,
\sigma\otimes\sigma\}$ has four braidings
$\e^+\otimes \e^+$, $\e^+\otimes \e^-$, 
$\e^-\otimes \e^+$ and $\e^-\otimes \e^-$.  The above 
consideration show that these four exhaust all the possibilities.

{\bf Subcase IIIc}: $D_{2n}$ with $n>3$.
We label the even vertices of $D_{2n}$ as
$0,2,4,\dots,2n-2_+,2n-2_-$ based on the $\alpha$-induction
as above.  The irreducible objects of the quantum double system 
are labeled with pairs of these again.
If $n>5$, then the statistical dimension of
$2$ is the smallest among $2,4,\dots,2n-2_+,2n-2_-$ 
and this shows that only possible embeddings arise
from $2\mapsto (0,2)$ and $2\mapsto (2,0)$.  For the
cases $n=4,5$, the fusion rules directly show that
these are also only possibilities.  They give 2 embeddings,
and the number of the braidings is again 2.

{\bf Case IV}: $E_6$.  It was found by Ocneanu that the system
arising from the even vertices of $E_6$ has no braiding.
This seems to be well-known to experts, but
we present arguments here as follows for the sake of completeness.

The dual principal graph of the Longo-Rehren subfactor
arising from the even vertices of $E_6$ is given in
\cite[Figure 1]{BEK3}.  With the labeling used in this
Figure, an embedding has to map $2$ to one of
$(2,0)$, $(8,0)$, $(1,1)$, $(5,1)_1$, $(5,1)_2$, but
none of these give the correct $E_6$ fusion rules.  For
example, we have $(2,0)\cdot(2,0)=(0,0)\oplus(2,0)\oplus(4,0)$.
This shows we have no embeddings.

{\bf Case V}: $E_8$.  This is very similar to Case IV and
the same remark applies.
We now use \cite[Figure 2]{BEK3}.  Any embedding has to
map $2$ to $(2,0)$, but this does not give the correct
$E_8$ fusion rules, so we have no embeddings.

Combining all the above, we obtain the following theorem.

\begin{theorem}
\label{num-braid}
For each system of bimodules given by the even vertices of
one of the $A_n$-$D_{2n}$-$E_{6,8}$ diagrams, the number of
braidings is given as in Table \ref{num-braid-tab}.

\begin{table}
\begin{center}
\begin{tabular}{|c|c|c|c|c|c|c|c|}\hline
$A_2$ & $A_5$ & other $A_n$ & $D_4$ & $D_6$ &
other $D_{2n}$ & $E_6$ & $E_8$\\ \hline
1 & 3 & 2 & 3 & 4 & 2 & 0 & 0 \\
\hline
\end{tabular}
\caption{The numbers of braidings}
\label{num-braid-tab}
\end{center}
\end{table}
\end{theorem}

\section{Another proof on the restriction on the index values}

We now give an alternative proof of the restriction on the
index values for inclusions of local conformal nets based
on the results in the previous section.

Suppose that one of the $A_n$-$D_{2n}$-$E_{6,8}$ diagrams
is realized as a principal graph
for an inclusion $\A(I)\subset\B(I)$ of
local conformal nets.  The even vertices of the principal
graph must give DHR-sectors of the net $\A$, so they
must have a braiding.  By Theorem \ref{num-braid},
the graphs $E_6$ and $E_8$ are excluded.

Suppose the graph is $A_n$, $n\neq 3, 5$.  By Theorem
\ref{num-braid}, all the braidings arise from $SU(2)_{n-1}$,
and $\A(I)\subset \B(I)$ is a local extension, so we can
copy the $Q$-system \cite{L}, \cite{LR} of $\A(I)$ so
that the net $SU(2)_{n-1}$ also has a local extension
with the dual canonical endomorphism $\theta=0\oplus 2$.
The classification table in \cite[Theorem 2.4]{KL} shows
that this is impossible.  So only $A_3$, $A_5$ remain,
and we already know they are indeed realized as in
Section \ref{alpha}.

Now suppose the graph is $D_{2n}$, $n\neq 2,3$.
Again by Theorem \ref{num-braid}, all the braidings
arise from the local extension of $SU(2)_{4n-4}$ with index 2,
and $\A(I)\subset \B(I)$ is a local extension, so we can
copy the $Q$-system of $\A(I)$ so
that the local extension of $SU(2)_{4n-4}$ also has a 
further local extension.  This is a local extension
of $SU(2)_{4n-4}$ with the dual canonical endomorphism
$\theta=0\oplus 2\oplus (4n-6)\oplus (4n-4)$.
The classification table in \cite[Theorem 2.4]{KL} shows
that this is again impossible.  So only $D_4$, $D_6$ remain,
and we already know they are indeed realized as in
Section \ref{alpha}.  Note that in the case of $D_6$,
the braiding arising from the local extension of $SU(2)_8$
with index 2 does not allow a further local extension,
but the representation category of $SU(2)_3\otimes
{\mathrm{Vir}}_{7/10}$ contains a subcategory corresponding
to the even vertices of $D_6$, and the braiding arising
from this local conformal net does allow a local extension
to $SU(2)_2\otimes SU(2)_1$.

We thus obtain the following theorem again.

\begin{theorem}
\label{rest-ind2}
The smallest four values of the Jones indices of inclusions
of local conformal nets are $1$, $2$, $3$,
$4\cos^2 \pi/10$ and $4$.
\end{theorem}
\medskip

\begin{remark}
{\rm In the above arguments on the restriction of possible index
values, all we need is a local $Q$-system arising from a unitary
braided tensor category.  So we have the same conclusion for
index values for such $Q$-systems.
}\end{remark}

\begin{remark}
{\rm We have a similar problem for index values also for inclusions
of 2-dimensional conformal nets.  The same arguments for the
restriction works.  The integer index values are again all possible
with orbifold nets.  

We see that the remaining value $4\cos^2 \pi/10$ is also possible
as follows.  Consider the $SU(2)_3$-net and label the irreducible
DHR sectors as $\lambda_0$, $\lambda_1$, $\lambda_2$, $\lambda_3$,
where $\lambda_0$ is the vacuum representation.  Then the Longo-Rehren
inclusion is given by a $Q$-system 
$\bigoplus_{j=0,1,2,3}\lambda_j \otimes \lambda^{\mathrm{opp}}_j$,
but the endomorphism 
$\bigoplus_{j=0,2}\lambda_j \otimes \lambda^{\mathrm{opp}}_j$ 
also gives a $Q$-system with localiy, hence a local extension.
This gives the index value $4\cos^2 \pi/10$.  So in this
setting, we have the same result as in
Theorem \ref{rest-ind2}.
}\end{remark}

\noindent
{\bf Acknowledgmentes.} This work has been partly done during visits of S.C. and R.L. at the University of Tokyo and of Y.K. at the University of Rome ``Tor Vergata''; we thank both institutions for the warm hospitality.


\begin{thebibliography}{99}

\bibitem{AH} M. Asaeda \& U. Haagerup, 
{Exotic subfactors of finite depth with Jones indices
${(5+\sqrt{13})}/{2}$ and ${(5+\sqrt{17})}/{2}$},
Commun. Math. Phys. {\bf 202} (1999) 1--63.

\bibitem{AY} M. Asaeda \& S. Yasuda,
{On Haagerup's list of potential principal graphs of subfactors},
Commun. Math. Phys. {\bf 286} (2009) 1141--1157.

\bibitem{BMPS} S. Bigelow, S. Morrison, E. Peters \& N. Snyder,
{Constructing the extended Haagerup planar algebra},
arXiv:0909.4099.

\bibitem{BE}
J. B\"ockenhauer \& D. E. Evans,
{Modular invariants, graphs and $\alpha$-induction for
nets of subfactors I},
Commun. Math. Phys. {\bf 197} (1998) 361--386. II
{\bf 200} (1999) 57--103. III {\bf 205} (1999) 183--228.

\bibitem{BEK1}
J. B\"ockenhauer, D. E. Evans \& Y. Kawahigashi,
{On $\alpha$-induction, chiral projectors
and modular invariants for subfactors},
Commun. Math. Phys. {\bf 208} (1999) 429--487.

\bibitem{BEK2}
J. B\"ockenhauer, D. E. Evans \& Y. Kawahigashi,
{Chiral structure of modular
invariants for subfactors}, 
Commun. Math. Phys. {\bf 210} (2000) 733--784.

\bibitem{BEK3}
J. B\"ockenhauer, D. E. Evans \& Y. Kawahigashi,
{Longo-Rehren subfactors arising from $\alpha$-induction},
Publ. RIMS Kyoto Univ. {\bf 37} (2001) 1--35.

\bibitem{DLR} C. D'Antoni, R. Longo \& F. R\u{a}dulescu, {Conformal nets, maximal 
temperature and models from free probability}, J. Operator Theory {\bf 45} (2001)
195--208. 

\bibitem{D} A. Degan,
``On the range of Jones index for inclusions of local nets of von Neumann 
algebras on the circle'',
Tesi di Dottorato, Universit\`a di Chieti-Pescara ``G. d'Annunzio'' 2009.

\bibitem{EK1}
D. E. Evans, Y. Kawahigashi, 
{Orbifold subfactors from Hecke algebras II
---Quantum doubles and braiding---},
Commun. Math. Phys. {\bf 196} (1998) 331--361

\bibitem{EK2} D. E. Evans \& Y. Kawahigashi,
``Quantum symmetries on operator algebras'',
Oxford University Press, 1998.

\bibitem{GHJ}
F. Goodman, P. de la Harpe \& V. F. R. Jones,
``Coxeter Graphs and Towers of Algebras''
MSRI Publications, Springer, Berlin-Heidelberg-New York,
{\bf 14}, 1989.

\bibitem{H}
U. Haagerup, 
{Principal graphs of subfactors in the index range 
$4< 3+\sqrt2$}, in {\em Subfactors ---
Proceedings of the Taniguchi Symposium, Katata ---},
(ed. H. Araki, et al.),
World Scientific (1994), 1--38.

\bibitem{I}
M. Izumi, 
{The structure of sectors associated with the Longo-Rehren
inclusions I. General theory}, 
Commun. Math. Phys. {\bf 213} (2000) 127--179. 

\bibitem{J}
V. F. R. Jones, {Index for subfactors}, Invent. Math. {\bf 72}
(1983) 1--25.

\bibitem{KL}
Y. Kawahigashi \& R. Longo,
{Classification of local conformal nets. Case $c<1$},
Ann. of Math. {\bf 160} (2004) 493--522.

\bibitem{KW} D. Kazhdan \& H. Wenzl,
{Reconstructing monoidal categories}.
In I. M. Gelfand Seminar, volume 16 of Adv.
Soviet Math., pages 111--136.
Amer. Math. Soc., Providence, RI, 1993.

\bibitem{L0}
R. Longo, {Minimal index and braided subfactors},
J. Funct. Anal. {\bf 109} (1992), 98-112.

\bibitem{L}
R. Longo, {A duality for Hopf algebras and for subfactors},
Commun. Math. Phys. {\bf 159} (1994) 133--150.

\bibitem{LR}
R. Longo \& K.-H. Rehren, {Nets of subfactors},
Rev. Math. Phys. {\bf 7} (1995) 567--597.

\bibitem{O}
A. Ocneanu,
{Chirality for operator algebras}
(Notes recorded by Y. Kawahigashi),
in {\em Subfactors}, World Scientific, (1994) 39--63

\bibitem{P0} S. Popa,
{Markov traces on universal Jones algebras and subfactors
of finite index},
Invent. Math. {\bf 111} (1993), 375--405. 

\bibitem{P} S. Popa,
``Classification of Subfactors and Their Endomorphisms'',
CBMS Regional Conference Series in Mathematics 1995.

\bibitem{R} K.-H. Rehren, 
Subfactors and coset models,
in {\em Generalized symmetries in physics}
(Clausthal, 1993),  World Scientific, (1994) 338--356.

\bibitem{R1}
K.-H. Rehren, {On the range of the index of subfactors},
J. Funct. Anal. {\bf 134} (1995) 183--193

\bibitem{X1}
F. Xu, 
{New braided endomorphisms from conformal inclusions},
Commun. Math. Phys. {\bf 192} (1998) 347--403.

\bibitem{X2}
F. Xu,
{Algebraic coset conformal field theories I},
Commun. Math. Phys. {\bf 211} (2000) 1--44.

\end{thebibliography}
\end{document}